\documentclass[a4paper,12pt]{amsart}
\usepackage{amssymb,amsthm}

\input xy
\xyoption{all}

\begin{document}

\newtheoremstyle{Def}{8pt}{8pt}{}{}{\bfseries}{.}{.5em}{\thmname{#1}\thmnote{ #3}}

\newcommand{\D}{\ensuremath{\mathcal{D}}}
\newcommand{\C}{\mathbb{C}}
\renewcommand{\P}{\mathbb{P}}
\newcommand{\SL}{\operatorname{SL}}
\renewcommand{\H}{\operatorname{H}}
\newcommand{\SHom}{\mathcal{H}om}
\newcommand{\F}{\ensuremath{\mathcal{F}}}
\newcommand{\SH}{\mathcal{H}}
\newcommand{\LW}{\mathcal{L}(v)}
\newcommand{\g}{\mathfrak{g}}
\newcommand{\height}{\operatorname{ht}}
\newcommand{\kar}{\operatorname{char}}
\newcommand{\Spec}{\operatorname{Spec}}
\newcommand{\Supp}{\operatorname{Supp}}
\newcommand{\Sing}{\operatorname{Sing}}
\newcommand{\Gr}{\operatorname{Gr}}
\newcommand{\cd}{\operatorname{cd}}

\newcommand{\Fbb}{\mathbb{F}}
\newcommand{\Gbb}{\mathbb{G}}
\newcommand{\Qbb}{\mathbb{Q}}
\newcommand{\Rbb}{\mathbb{R}}
\newcommand{\Zbb}{\mathbb{Z}}
\newcommand{\Lcal}{\mathcal L}
\newcommand{\Mcal}{\mathcal M}
\newcommand{\Ocal}{\mathcal O}
\newcommand{\diag}{\operatorname{diag}}
\newcommand{\GL}{\operatorname{GL}}

\newtheorem{Corollary}{Corollary}[section]
\newtheorem{Lemma}{Lemma}[section]
\newtheorem{Theorem}{Theorem}[section]
\newtheorem{Proposition}[Theorem]{Proposition}
{\theoremstyle{Def} \newtheorem{Definition}{Definition}}
\newtheorem{Remark}{Remark}
\newenvironment{Proof}{{\it Proof.\/}}{\hfill $\square$\medskip}

\title[Frobenius splitting of regular conjugacy classes]
{Frobenius splitting of equivariant closures of regular
conjugacy classes}
\author{Jesper Funch Thomsen}
\address{Institut for matematiske fag\\Aarhus Universitet\\ 8000 \AA rhus C, Denmark}
\email{funch@imf.au.dk}

\begin{abstract}
Let $G$ denote a connected semisimple and simply 
connected algebraic group over an algebraically 
closed field $k$ of positive characteristic 
and let $g$ denote a regular element of $G$. 
Let $X$ denote any equivariant embedding of $G$. 
We prove that the closure of the conjugacy 
class of $g$ within $X$ is normal and 
Cohen-Macaulay. Moreover, when $X$ is smooth 
we prove that this closure is a local complete
intersection. As a consequence, the closure of 
the unipotent variety within $X$ share the
same geometric properties.
\end{abstract}

\maketitle

\section{Introduction}

Let $G$ denote a connected semisimple  
linear algebraic group over 
an algebraically closed field $k$ of positive
characteristic. Consider 
$G$ as a $G \times G$-variety by left and 
right translation. An equivariant $G$-embedding 
is a normal $G \times G$-variety $X$ containing 
a dense open subset $G \times G$-equivariantly 
isomorphic to $G$. Let $B$ denote
a Borel subgroup of $G$ and let $g$ denote 
an element of $G$. The closure of the double 
coset $B g B$ within $X$ is called a {\it large 
Schubert variety}. As proved in \cite{BriTho}
(see also \cite{BriPol}) large Schubert varieties 
are normal and Cohen-Macaulay. In the present 
paper we will prove a similar result for 
closures of diagonal $G$-orbits of 
regular elements within $X$.

Assume that $G$ is simply connected. An element
$g$ in $G$ is called regular if the centralizer
of $g$ in $G$ is of dimension equal to the rank
$l$ of $G$. By \cite{Steinberg} the set of 
$G$-conjugacy classes of regular elements is 
parametrized by the set of points in affine 
space $\mathbb A^l$. To describe this 
correspondence, let $\chi_1, \dots, \chi_l$
denote the set of fundamental $G$-characters 
and consider the Steinberg map $\chi :
G \rightarrow \mathbb A^l$ defined by 
$\chi(h) = (\chi_1(h), \dots, \chi_l(h))$. 
Then the fiber $\chi^{-1}({\bf a})$, 
associated to a point ${\bf a}$ in $\mathbb A^l$, 
contains a unique open dense $G$-orbit 
consisting of the set of regular elements 
within $\chi^{-1}({\bf a})$. Hence,  
$\chi^{-1}({\bf a})$ is the closure of
a unique regular conjugacy class within $G$.

Let $g$ denote any regular element within $G$.
We prove that the closure $\overline{C_G(g)}$ 
of the $G$-conjugacy class of $g$ within an 
equivariant embedding $X$ of $G$ is normal and 
Cohen-Macaulay. Moreover, when $X$ is smooth 
we prove that $\overline{C_G(g)}$ is a local 
complete intersection and we calculate
its associated dualizing sheaf. When $X=G$ 
this statement is equivalent to saying that 
the fibers of the Steinberg map are all normal, 
Cohen-Macaulay and local complete intersections. 
The latter statement is due to Steinberg 
(Thm.6.11 and Thm.8.1 in \cite{Steinberg}).
The case of primary interest to us is when
$g$ is a regular unipotent element in $G$. 
In this case  $\overline{C_G(g)}$ coincides 
with the closure of the unipotent variety 
$\mathcal U$ within $X$. In particular, when
$G$ coincides with the associated group of
adjoint type (i.e. type $E_8$, $F_4$ and 
$G_2$) we obtain a description of the 
geometry of the closure of the unipotent 
variety within the wonderful compactification 
of $G$.

In order to prove the described results we
will use the theory of Frobenius splitting. 
It is well known that any equivariant embedding
$X$ of $G$ is Frobenius split. We prove that
there exists a Frobenius splittings of $X$
which compatibly splits the closure 
$\overline{C_G(g)}$ of a given regular conjugacy 
class. Moreover, we may choose such a 
splitting to be canonical in the sense 
of Mathieu (see Chap. 4, \cite{BriKum}). In particular, when 
$\mathcal L$ is a $G$-linearized line bundle on 
$\overline{C_G(g)}$, then the associated set of global 
sections ${\rm H}^0(\overline{C_G(g)}, \mathcal L)$
admits a good filtration (as a $G$-module).
Finally, we prove that when $X$ (and hence 
$\overline{C_G(g)}$) is projective then the
higher cohomology groups of globally generated 
line bundles on $\overline{C_G(g)}$ are zero.

\section{Notation}
\label{notation}

In the following sections $G$ will denote a connected
semisimple linear algebra group over an algebraically
closed field $k$ of positive characteristic $p>0$.
The associated group of adjoint type will be denoted 
$G_{\rm ad}$. We will fix a Borel subgroup $B$ and a
maximal torus $T$ within $B$. The character (resp. 
co-characters) of $T$ will be denoted by $X^*(T)$
(resp. $X_*(T)$). We identify $X^*(T)$ with the 
group of $B$-characters $X^*(B)$.

The set of roots associated with $T$ will be denoted 
by $R$. We define a root to be positive if the associated 
$T$-weight space (under the adjoint action) in the Lie 
algebra of $B$ is zero. The set of positive roots 
is denoted by $R^+$ and the associated set of simple 
roots is denoted by $\Delta = \{ \alpha_1, \dots, 
\alpha_l \}$. To each root $\alpha$ in $R$ there is 
an associated coroot $\alpha^\vee$ which we use 
to define the set of dominant weights in $X^*(T)$. 
When $G$ is simply connected every dominant weight
is a positive linear combination of the fundamental 
dominant weights $\omega_1, \dots, \omega_l$, 
where $\omega_i$ denotes  the fundamental weight 
associated to $\alpha_i$.

The Weyl group $W$ associated with $T$ parametrizes 
the set of Schubert varieties $X(w)$ in the flag 
variety $G/B$. When $w \in W$ is represented by 
an element $\dot{w} \in G$ we will write $B w B$ 
for the subset $B \dot{w} B$ of $G$. The simple 
reflections, denoted by $s_1, \dots, s_l$, generate 
the Weyl group and we use them to define the length $l(w)$
of an element $w \in W$. The unique element in $W$ of 
maximal length is denoted by $w_0$.

For any $B$-character $\lambda$ 
we define a $G$-module by
$${\rm H}(\lambda) = \{ f : G \rightarrow k : 
f(gb) = \lambda(b)^{-1} f(g) , ~b \in B, g \in G
\}.$$
Then ${\rm H}(\lambda)$ is nonzero exactly when 
$\lambda$ is a dominant $T$-character. When 
${\rm H}(\lambda)$ is nonzero there exists a 
unique $B$-stable line within ${\rm H}(\lambda)$. 
This line is called a lowest weight line and 
its associated weight is $w_0 \lambda$. Similarly, 
${\rm H}(\lambda)$ contains a unique highest
weight line invariant under the opposite Borel 
group $B^+$ and of weight $\lambda$.
The $G$-module ${\rm H}(\lambda)$ coincides 
with the global sections of a unique $G$-linearized
line bundle on $G/B$ which we denote by 
${\mathcal L}_{G/B}(\lambda)$.

When $M$ is a $G$-module we consider the group 
of endomorphisms ${\rm End}(M)$ as a 
$G \times G$-module by 
$$( (g',g) \cdot f )(m) = g' \cdot (f(g^{-1} \cdot m))$$
for $g, g' \in G, f \in {\rm End}(M)$ and $m \in M$. 
When considering  ${\rm End}(M)$ as a $G$-module
we do this by identifying $G$ with the subgroup
 $G \times \{ e \} $ of $G \times G$. The field
$k$ will always  be considered as the trivial 
representation.

\subsection{The Steinberg map}

Let now $G$ denote a simply connected group and let 
$M$ denote a finite dimensional rational $G$-module
defined by a morphism : $G \rightarrow {\rm GL(M)}$. 
The $G$-character $\chi_M :G \rightarrow k$ of $M$ 
is by definition the composition of latter map with the
trace function on  ${\rm GL(M)}$. The 
{\it Steinberg map} (\cite{Steinberg}) is the map
$$ \chi : G \rightarrow {\mathbb A}^l,$$
$$ g \mapsto (\chi_1(g), \dots, \chi_l(g)),$$
defined as the product of the $G$-characters
$\chi_i$ associated to the $G$-modules 
${\rm H}(\omega_i)$. When ${\bf a} = (a_1,
\dots, a_l) \in \mathbb A^l$ we will by 
$\chi^{-1}({\bf a})$ denote the (scheme theoretic) 
fiber of $\chi$ at the point ${\bf a}$. As
mentioned in the introduction we may consider 
the fibers of the Steinberg map as closures 
of conjugacy classes of regular elements within 
$G$. Among all fibers the fiber at ${\bf a} = 
(\chi_1(e), \dots, \chi_l(e))$ is of particular 
importance as this fiber coincides with the set 
of unipotent elements within $G$.

\section{Equivariant embeddings}

In this section $G$ denotes a connected semisimple linear
algebraic group. We think of $G$ as a $G \times G$-variety
by left and right translation. An {\it equivariant 
$G$-embedding} (or simply a $G$-embedding) is a normal 
$G \times G$-variety $X$ containing an open subset 
which is $G \times G$-equivariantly isomorphic to $G$.

\subsection{The wonderful compactification}

When $G=G_{\rm ad}$ is of adjoint type there exists a 
distinguished equivariant embedding ${\bf X}$ of 
$G$ which is called the {\it wonderful compactification} 
(see e.g. Section 6.1,\cite{BriKum}).

The wonderful compactification ${\bf X}$ is
a smooth projective variety such that the
complement of the open subset $G$ is a finite
union of smooth irreducible divisors which 
intersect transversally. Moreover, the 
common intersection of the irreducible components 
of ${\bf X} \setminus G$ is a closed 
$G \times G$-orbit ${\bf Y}$ which is 
$G \times G$-equivariantly isomorphic to 
$G/B \times G/B$. The variety ${\bf Y}$ 
is the unique closed $G \times G$-orbit 
within $X$.

\subsection{Toroidal embeddings}
\label{toi-emb} 
An embedding $X$ of a semisimple
group $G$ is called {\it toroidal} if the 
canonical map $\phi : G \rightarrow G_{\rm ad}$ 
admits an extension 
$\pi_X : X \rightarrow  {\bf X}$
into the  wonderful compactification ${\bf X}$
of the group  $G_{\rm ad}$ of adjoint
type. 

When $X$ is a complete toroidal embedding 
then every closed orbit of $X$ will
map to the unique closed orbit ${\bf Y} 
\simeq G/B \times G/B$ within ${\bf X}$. 
As a consequence, every closed orbit of 
$X$ is then $G \times G$-equivariantly 
isomorphic to $G/B \times G/B$.

A toroidal embedding is uniquely determined
by the closure $\overline T$ of $T$ within 
$X$. The closure $\overline T$ is normal 
and hence $\overline T$ is a toric variety 
with respect to $T$. As toric varieties admits 
resolutions
by smooth toric varieties this essentially
explains the following (see also Prop.6.2.5,\cite{BriKum})

\begin{Theorem}
\label{resolution}
For any equivariant $G$-embedding $X$
there exists a smooth toroidal embedding 
$X'$ of $G$ and a birational projective
morphism $X' \rightarrow X$ extending 
the identity map on $G$.
\end{Theorem}

\section{Frobenius splitting}

Let $X$ denote a scheme of finite type over 
an algebraically closed field $k$ of 
characteristic $p>0$. The {\it absolute Frobenius morphism} 
on $X$ is the morphism $F : X \rightarrow X$ of schemes, 
which is the identity on the set of points and where 
the associated map of sheaves 
$$ F^\sharp : \mathcal O_X  \rightarrow 
F_* \mathcal O_X$$
is the $p$-th power map. We say that $X$ is
{\it Frobenius split} if there exists a morphism
$s \in {\rm Hom}_{\mathcal O_X} (F_* \mathcal O_X,  
\mathcal O_X )$  such that the composition 
 $s \circ  F^\sharp$ is the identity 
map on $\mathcal O_X$.

\subsection{Stable Frobenius splittings along divisors}

Let $D$ denote an effective Cartier divisor on $X$
with associated line bundle $\mathcal O_X(D)$ and 
canonical section $\sigma_D$. We say that $X$ is 
{\it stably Frobenius split along $D$} if there 
exists a positive integer $e$ and an morphism 
$$ s \in {\rm Hom}_{\mathcal O_X} (F_*^e \mathcal O_X(D),
\mathcal O_X ),$$
such that $s(\sigma_D)=1$. In this 
case we say that $s$ is a stable Frobenius splitting
of $X$ along $D$ of degree $e$. Notice that $X$ is 
Frobenius split exactly when there exists a stable 
Frobenius splitting of $X$ along the zero divisor $D=0$.

\begin{Remark}
\label{opensubset}
Consider an element $s \in {\rm Hom}_{\mathcal O_X} 
(F_*^e \mathcal O_X(D), \mathcal O_X )$. Then the 
condition $s(\sigma_D)=1$ on $s$ for it to define 
a stable Frobenius splitting of $X$, may be checked 
on any open dense subset of $X$.  
\end{Remark}

\subsection{Subdivisors}
\label{subdivisor}

Let $D' \leq D$ denote an effective Cartier 
subdivisor and let $s$ be a stable Frobenius 
splitting of $X$ along $D$ of degree $e$. 
The composition of $s$ with the map
$$ F_*^e \mathcal O_X(D') \rightarrow 
F_*^e \mathcal O_X(D),$$
defined by the canonical section of the divisor
$D-D'$, is then a stable Frobenius splitting of
$X$ along $D'$ of degree $e$. Applying this 
to the case $D'=0$ it follows that if $X$ is 
stably Frobenius split along any effective divisor 
$D$ then $X$ is also Frobenius split.

\begin{Lemma}
\label{Frob1}
Let $D_1$ and $D_2$ denote effective Cartier divisors.
If $s_1$ (resp. $s_2$) is a stable Frobenius splitting
of $X$ along $D_1$ (resp. $D_2$) of degree $e_1$ (resp.
$e_2$), then there exists a stable Frobenius splitting 
of $X$ along $D_1 + D_2$ of degree $e_1 + e_2$.
\end{Lemma}
\begin{proof}
By the discussion above it suffices to prove that 
there exists a stable Frobenius splitting of $X$ along
$D_1 + p^{e_1} D_2$ of degree $e_1 + e_2$. By the 
projection formula there exists an identification 
$$F_*^{e_1 + e_2}\mathcal O_X(D_1 + p^{e_1} D_2)
\simeq F_*^{e_2} ( \mathcal O_X(D_2) \otimes 
F_*^{e_1} \mathcal O_X(D_1) ).$$ 
where we have used the relation $F^* O_X(D_2)
\simeq O_X(p D_2)$. Hence, $s_1$ defines a map 
$$F_*^{e_1 + e_2}\mathcal O_X(D_1 + p^{e_1} D_2)
\rightarrow F_*^{e_2} ( \mathcal O_X(D_2)).$$
Composing the latter map with $s_2$ we 
obtain a map 
$$F_*^{e_1 + e_2}\mathcal O_X(D_1 + p^{e_1} D_2)
\rightarrow \mathcal O_X,$$
which defines the desired stable Frobenius splitting 
of $X$.
\end{proof}

\subsection{Compatibly split subschemes}

Let $Y$ denote a closed subscheme of $X$ with
sheaf of ideals $\mathcal I_Y$. When 
$$s \in {\rm Hom}_{\mathcal O_X} (F_*^e \mathcal O_X(D),
\mathcal O_X ),$$
is a stable Frobenius splitting of $X$ along $D$ we 
say that $s$ {\it compatibly Frobenius splits} 
$Y$ if the following conditions are satisfied

\begin{enumerate}
\item The support of $D$ does not contain any
of the irreducible components of $Y$.
\item $s(F_*^e( \mathcal I_Y \otimes \mathcal O_X(D)))
\subseteq \mathcal I_Y$.
\end{enumerate}
When $s$ compatibly Frobenius splits $Y$
there exists an induced stable Frobenius 
splitting of $Y$ along $D \cap Y$ of degree
$e$. Notice that when only condition (2)
is satisfied then $Y$ is still compatibly 
Frobenius split by the induced stable 
Frobenius splitting of $X$ along the zero 
divisor $0$. In concrete situations 
condition (2) may be checked using the following 
result

\begin{Lemma}
\label{Frob1.5}
Let $s$ denote a stable Frobenius splitting of 
$X$ along a divisor $D$ and let $Y$ denote a 
closed subscheme of $X$ satisfying the above 
condition (1). If $Y$ is compatibly Frobenius
split by the induced stable Frobenius splitting 
of $X$ along the zero divisor $0 \leq D$,
then $Y$ is also compatibly Frobenius split 
by $s$.
\end{Lemma}
\begin{proof}
Argue as in the proof of Prop.1.4 in \cite{Ram}.   
\end{proof}

\begin{Lemma}
\label{Frob2}
Let $s$ denote a stable Frobenius splitting of $X$ 
along $D$ which compatibly Frobenius
splits a closed subscheme $Y$ of $X$. If $D' \leq
D$ then the induced stable Frobenius splitting of $X$ 
along $D'$, defined in Section \ref{subdivisor},     
compatibly Frobenius splits $Y$. 
\end{Lemma}
\begin{proof}
This follows immediately from the construction of 
the induced Frobenius splitting of $X$ along $D'$. 
\end{proof}

\begin{Lemma}
\label{Frob3}
Let $D_1$ and $D_2$ denote effective Cartier divisors.
If $s_1$ (resp. $s_2$) is a stable Frobenius splitting
of $X$ along $D_1$ (resp. $D_2$) of degree $e_1$ (resp.
$e_2$) which compatibly splits a closed subscheme $Y$
of $X$, then there exists a stable Frobenius splitting 
of $X$ along $D_1 + D_2$ of degree $e_1 + e_2$ which 
compatibly splits $Y$.
\end{Lemma}
\begin{proof}
The stable Frobenius splitting of $X$ along $D_1+D_2$
defined in the proof of Lemma \ref{Frob1} compatibly
Frobenius splits $Y$.
\end{proof}

As an easy consequence of the above definitions we find

\begin{Lemma}
\label{Frob4}
Let $s$ denote a stable Frobenius splitting of 
$X$ along an effective divisor $D$. Then
\begin{enumerate}
\item $X$ is reduced and every irreducible component 
of $X$ is compatibly Frobenius split. 
\item If $s$ compatibly Frobenius splits a closed 
subscheme $Y$ of $X$ then each irreducible component of
$Y$ is also compatibly Frobenius split by $s$.
\item Assume that $s$ compatibly Frobenius splits closed
subschemes $Y_1$ and $Y_2$ and that the support of $D$
does not contain any of the irreducible components of 
the scheme theoretic intersection $Y_1 \cap Y_2$. Then 
$s$ compatibly Frobenius splits  $Y_1 \cap Y_2$. 
\end{enumerate}
\end{Lemma}

The following statement relates stable Frobenius 
splitting along divisors with compatibly Frobenius
splitting.

\begin{Lemma}
\label{Frob5}
Let $D$ and $D'$ denote effective Cartier divisors 
and let $s$ denote a stable Frobenius splitting of 
$X$ along $(p-1)D + D'$ of degree $1$. Then there
exists a stable Frobenius splitting of $X$ along 
$D'$ of degree $1$ which compatibly splits the closed
subscheme defined by $D$. 
\end{Lemma}
\begin{proof}
Let $s$ denote a stable Frobenius splitting
along $(p-1)D + D'$ of degree $1$. Define $s'$
to be the composition of $s$ with the map 
$$ F_* \mathcal O_X(D') \rightarrow F_* 
\mathcal O_X((p-1)D +D'),$$
defined by the canonical section $\sigma_{D}^{p-1}$ 
of $(p-1)D$.  Then $s'$ is a stable 
Frobenius splitting of $X$ of degree $1$. 
It remains to show that $D$ is compatibly 
Frobenius split by $s'$. As this is a local condition
we may assume that $\mathcal O_X(D)$ and 
$\mathcal O_X(D')$ are trivial line bundles and 
that $X$ is affine. Identify $\sigma_D$ and 
$\sigma_{D'}$ with elements in the coordinate
ring $k[X]$ and consider $s$ and $s'$ as
maps from $F_ * k[X]$ to $k[X]$. By definition
$$ s'(a) = s( a \sigma_D^{p-1})~,~a \in k[X].$$
Hence, when $a = \sigma_D a' \in (\sigma_D)$ belongs 
to the ideal generated by $\sigma_D$ we find 
$$s'(a) = s(\sigma_D^p a') = \sigma_D s(a') \in (\sigma_D).$$
It remains to prove that none of the components
of $D$ is contained in the support of $D'$.
Restricting, if necessary, to an open subset 
we may assume that $D$ is irreducible and 
nonempty. Now assume that $D$ is contained in the support 
of $D'$. By Lemma \ref{Frob4} the ideal 
$(\sigma_D)$ is radical. Hence, $\sigma_{D'}$ is
contained in $(\sigma_D)$. In particular, 
$$ 1 = s'(\sigma_{D'}) \in (\sigma_D),$$ 
which is a contradiction. 
\end{proof}

\subsection{Cohomology and Frobenius splitting}

The notion of Frobenius splitting is particular 
useful in connection with proving higher cohomology 
vanishing for line bundles. The main idea 
is that when $s \in {\rm Hom}_{\mathcal O_X} 
(F_*^e \mathcal O_X(D), \mathcal O_X )$ is 
a stable Frobenius splitting of $X$ along 
the divisor $D$, then $s$ defines a splitting 
of the injective map 
$$ \mathcal O_X \rightarrow F_*^e \mathcal O_X(D).$$
Tensoring the latter map with a line bundle 
$\mathcal L$ on $X$ we find a split injective map
$$ \mathcal L \rightarrow F_*^e (\mathcal L^{p^e} 
\otimes \mathcal O_X(D)),$$
where we have applied the projection formula and
the relation $F^* \mathcal L = \mathcal L^p$.
As $F$ is a finite morphism the following
statement is an easy consequence.

\begin{Lemma}
\label{Frob6}
Let $s$ denote a stable Frobenius splitting 
of $X$ along $D$ of degree $e$. Then for every
line bundle $\mathcal L$ on $X$ and 
every integer $i$ there exists an 
inclusion 
$${\rm H}^i(X, \mathcal L) \subseteq 
{\rm H}^i(X, \mathcal L^{p^e} \otimes 
\mathcal{O}_X(D)),$$
of abelian groups. In particular, when
$X$ is projective, $\mathcal L$ is 
globally generated and $D$ is ample 
then the group ${\rm H}^i(X, \mathcal L)$ 
is zero for $i>0$.
\end{Lemma}

Assume now, moreover, that $Y$ is compatibly 
Frobenius split by the stable Frobenius splitting 
$s$ of $X$. Then $s$ defines a splitting of the map 
$$ \mathcal I_Y \rightarrow 
 F_*^e (\mathcal I_Y  \otimes \mathcal O_X(D)),$$
$$ f \mapsto f^{p^e} \sigma_D.$$
Tensoring the latter map with a line bundle
$\mathcal L$ on $X$ then leads to 

\begin{Lemma}
\label{Frob7}
Let $s$ denote a stable Frobenius splitting 
of $X$ along $D$ of degree $e$  and let $Y$
denote a closed compatibly Frobenius split 
subscheme of $X$. Then for every
line bundle $\mathcal L$ on $X$ and every 
integer $i$ there exists an inclusion
$${\rm H}^i(X, \mathcal I_Y \otimes 
\mathcal L) \subseteq 
{\rm H}^i(X, \mathcal I_Y \otimes 
\mathcal L^{p^e} \otimes \mathcal{O}_X(D)),$$
of abelian groups. In particular, when
$X$ is projective, $\mathcal L$ is 
globally generated and $D$ is ample 
then the group ${\rm H}^i(X, \mathcal I_Y \otimes 
\mathcal L)$ is zero for $i>0$.
\end{Lemma}

\subsection{Push forward}
\label{push-forward}

Let $f : X \rightarrow X'$ denote proper 
morphism of schemes and assume that the 
induced map $ \mathcal O_{X'} \rightarrow 
f_* \mathcal O_X$ is an isomorphism. Then 
every Frobenius splitting of $X$ induces,
by application of the functor $f_*$, 
a Frobenius splitting of $X'$. Moreover, 
when $Y$ is a compatibly Frobenius 
split subscheme of $X$ then the induced 
Frobenius splitting of $X'$ compatibly
splits the scheme theoretic image
$f(Y)$ (Prop.4, \cite{MehtaRamanathan}).

We will need the following connected 
statement.

\begin{Lemma}
\label{Frob8}
Let $f : X \rightarrow X'$ denote a 
morphism of projective schemes such 
that $ \mathcal O_{X'} \rightarrow 
f_* \mathcal O_X$ is an isomorphism. 
Let $Y$ be a closed subscheme of 
$X$ and denote by $Y'$ the scheme
theoretic image $f(Y)$. Assume 
that there exists a stable Frobenius
splitting of $X$ along an ample 
divisor $D$ which compatibly splits 
$Y$. Then $f_* \mathcal O_Y 
= \mathcal O_{Y'}$ and $R^i f_* 
\mathcal O_Y = 0$ for $i>0$.
\end{Lemma}
\begin{proof}
Let $\mathcal L$ denote an ample line bundle on $X'$. 
By a result of G. Kempf (see Lemma 2.11, \cite{Ram})
it suffices to prove, for sufficiently large values 
of $n$, that
\begin{enumerate}
\item ${\rm H}^i(X, f^* \mathcal L^n) = 
{\rm H}^i(Y, f^* \mathcal L^n)= 0$ for $i>0$. 
\item The restriction map ${\rm H}^0(X, f^* \mathcal L^n) 
\rightarrow {\rm H}^0(Y, f^* \mathcal L^n)$ is surjective.
\end{enumerate}
Now apply the ``in particular'' statements in 
Lemma \ref{Frob6} and \ref{Frob7}.
\end{proof}

\subsection{Frobenius splitting of smooth varieties}

When $X$ is smooth there exists a 
canonical $\mathcal O_X$-linear identification 
(see e.g. Section 1.3.7 in \cite{BriKum}) 
$$ C_X : F_* \omega_X^{1-p} \simeq 
{\rm Hom}_{\mathcal O_X} (F_* \mathcal O_X,  
\mathcal O_X ).$$
Hence, a Frobenius splitting of $X$ may 
be identified with a global section of 
$\omega_X^{1-p}$ with certain properties.
A global section $s$ of $\omega_X^{1-p}$
which corresponds to a Frobenius splitting 
will be called a Frobenius splitting.

\begin{Lemma}
\label{Frob9}
Assume that $X$ is smooth and let $\tau$
denote a global section of $\omega_X^{1-p}$
which defines a Frobenius splitting of $X$.  
Then there exists a stable Frobenius 
splitting of $X$ of degree $1$ along 
the Cartier divisor defined by $\tau$. In 
particular, if $\tau = \tilde \tau^{p-1}$ 
is a $(p-1)$-th power of a global section 
$\tilde \tau $ of $\omega_X^{-1}$, then $X$ 
is Frobenius split compatibly with the zero 
divisor of $\tilde \tau$.
\end{Lemma}
\begin{proof}
Composing the evaluation map 
$$ {\rm Hom}_{\mathcal O_X} (F_* \mathcal O_X,  
\mathcal O_X ) \rightarrow \mathcal O_X,$$
$$ s \mapsto s(1),$$
with $C_X$ defines a stable Frobenius
splitting of $X$ of degree $1$ along 
the Cartier divisor defined by $\tau$. This 
proves the first assertion. The second 
assertion follows from Lemma \ref{Frob5}
with $D'=0$.
\end{proof}

\subsection{Frobenius splitting of $G/B$}

Let $G$ denote a simply connected linear algebraic
group. The flag variety $X=G/B$ is a smooth 
variety with dual canonical bundle 
$$\omega_X^{-1} = \mathcal L_{G/B}(2 \rho),$$
where $\rho$ denotes the sum of the fundamental 
dominant weights. Consider the multiplication
map 
$$ \phi_{G/B} : {\rm St} \otimes {\rm St} \rightarrow 
 {\rm H}(2(p-1)\rho) = {\rm H}^0(G/B, 
\omega_{G/B}^{1-p}),$$
where we have used the notation ${\rm St}$
to denote the Steinberg module ${\rm H}((p-1)\rho)$.
The Steinberg module ${\rm St}$ is an irreducible
selfdual $G$-module and hence there exists
a unique (up to nonzero scalars) nondegenerate 
$G$-invariant bilinear form $(,)$ on ${\rm St}$. Now we have 
the following description

\begin{Theorem} (\cite{LauTho})
\label{splitGB}
Let $h= \sum_i f_i \otimes g_i$ denote an element 
in ${\rm St} \otimes {\rm St}$. Then $\phi_{G/B}
(h)$ defines a Frobenius splitting of $X$ (up to a
nonzero scalar) if and only if $\sum_i (f_i,g_i)
\neq 0$.   
\end{Theorem}

As an easy consequence of this result (with $G$ 
substituted with $G \times G$) we have

\begin{Corollary}(\cite{LauTho})
\label{splitdiag}
Let $v_\Delta$ denote a generator of the unique 
diagonal $G$-invariant line in $ {\rm St} 
\otimes {\rm St}$ and let $v_-$ (resp. $v_+$) 
denote a generator of the highest (resp.
lowest) weight line in ${\rm St}$. Then
$$\phi_{(G \times G)/(B \times B)}(  v_\Delta 
\otimes (v_- \otimes v_+))$$
defines a Frobenius splitting of 
$G/B \times G/B$ (up to a nonzero 
scalar).
\end{Corollary}

\section{Preliminary results} 
\label{psi} 

Throughout this section $G$ is simply connected and
$\mu$ is a nonzero fixed dominant $T$-character.

\subsection{The morphism $\psi_{\mu}$} 
Consider the $G \times G$-equivariant defined by 
$$ \psi_{\mu} : G \rightarrow {\mathbb P}_\mu 
:={\mathbb P}({\rm End}({\rm H}(\mu))
\oplus k), $$
$$ g \mapsto [(g \cdot I_{{\rm H}(\mu)},1)],$$
where $I_{{\rm H}(\mu)}$ denotes the 
identity map on ${\rm H}(\mu)$. Let 
${\mathcal O}_\mu(1)$ denote the ample generator
of the Picard group of ${\mathbb P}_\mu$. The 
pull back $\psi_\mu^*({\mathcal O}_\mu(1))$ is
then (by the description of $\psi_\mu$) 
canonical isomorphic to the trivial line
bundle on $G$. The induced map on global 
sections  
$$ \psi_{\mu,G}^* : {\rm End}({\rm H}(\mu))^*
\oplus k \simeq {\rm H}^0( {\mathbb P}_\mu,
{\mathcal O}_\mu(1)) \rightarrow 
{\rm H}^0(G,
\psi_\mu^*({\mathcal O}_\mu(1))) \simeq k[G],$$
is given by 
$$\psi_ {\mu,G}^*(v^* \otimes v,a)(g) = v^*(g v) + a,$$
with $v^* \otimes v \in {\rm H}(\mu)^*
\otimes {\rm H}(\mu) \simeq {\rm End}({\rm H}(\mu))^*$
and $a \in k$.

Let $v_\mu^-$ (resp. $u_\mu^-$ ) denote a lowest 
weight vector in ${\rm H}(\mu)$ (resp. ${\rm H}(\mu)^*$).
Similarly we let  $v_\mu^+$ (resp. $u_\mu^+$ ) denote a 
highest weight vector in ${\rm H}(\mu)$ (resp. 
${\rm H}(\mu)^*$). Let ${\rm Tr}_\mu$ denote the 
trace function on ${\rm End}({\rm H}(\mu))$. Then

\begin{Lemma}
\label{OnG}
The function $\psi_{\mu,G}^*(u_\mu^- \otimes v_\mu^-,0)$
on $G$ is a generator of the lowest weight space in 
${\rm H}(- w_0\mu) \subseteq k[G]$. Moreover, the function
$\psi_{\mu,G}^*({\rm Tr}_\mu,0)$ coincides with the
$G$-character $\chi_\mu$ of ${\rm H}(\mu)$. 
\end{Lemma} 
\begin{proof}
That $\psi_{\mu,G}^*(u_\mu^- \otimes v_\mu^-,0)$ is nonzero
follows by evaluating it at $w_0$. The first 
assertion then follows as $u_\mu^- \otimes v_\mu^-$
is $B \times B$-semiinvariant (i.e. invariant up
to scalars) of weight $(-\mu, w_0 \mu)$.
The second assertion is clear from the  
discussion above.
\end{proof}

\subsubsection{Regular functions on $G$}
\label{functions}

In connection with Lemma \ref{OnG} we will 
later need a more precise description of 
$\psi_{\mu,G}^*(u_\mu^- \otimes v_\mu^-,0)$
when $\mu$ is a fundamental dominant weight.
Recall that the coordinate ring $k[G]$ is a 
unique factorization domain as $G$ is assumed 
to be simply connected. Then 

\begin{Lemma}
\label{lowestfunction}
Let $f_i$ denote a generator of the lowest weight 
space in ${\rm H}(\omega_i)$. Then $f_i$ is an 
irreducible element in $k[G]$ and the ideal 
$(f_i)$ generated by $f_i$ in $k[G]$ coincides 
with the ideal of functions vanishing on the 
closure $\overline{B w_0 s_i B}$  of the 
Bruhat cell $B w_0 s_iB$.
\end{Lemma}
\begin{proof}
Write $f_i$ as a product of irreducible 
elements
$f_i = f_{i,1} \cdots f_{i,n}.$
As $f_i$ is $B \times B$-semiinvariant each 
factor $f_{i,j}$ of this product 
is also $B \times B$-semiinvariant. 
Hence, there exists dominant weights 
$\omega_{i,j}$ such that $f_{i,j}$ 
generates the lowest weight space in 
${\rm H}({\omega_{i,j}}).$
Moreover, by the factorization of $f_i$
we must have $\omega_i = \sum_{j=1}^n 
\omega_{i,j}$. But then $n$ must 
be equal to 1 which proves the first
part of the statement.

As $f_i$ is irreducible the ideal $(f_i)$
is a prime ideal of height 1. By the 
$B \times B$-semiinvariance of $f_i$ this
ideal must coincide with the ideal of 
functions vanishing on the closure of some 
Bruhat cell $\overline{B w_0 s_j B}$
of codimension 1. It remains to prove 
that $f(w)= 0$, where $w \in G$ denotes 
a representative for $w_0 s_i$. So assume 
that $f(w) \neq 0$. When $t \in T$ we 
may calculate $ f(t w)$ in two different
ways. First by definition of ${\rm H}
(\omega_i)$ :
$$ f(t w) = f(w w^{-1} t w) = 
\omega_i(w^{-1} t^{-1} w) f(w) = 
(-w_0 s_i \omega_i)(t)f(w).$$
But also $f(t w) = (t^{-1}\cdot f)(w) = 
(-w_0 \omega_i)(t) f(w)$. In particular,
we conclude that $w_0 s_i \omega_i = 
w_0 \omega_i$ which is a contradiction.
\end{proof}

\subsection{The closed orbit in $Z_\mu$}
\label{Y}

In the sequel we let $Z_\mu$ denote
the closure of the image of $\psi_\mu$. 

\begin{Lemma}
\label{fixed}
The element $[(v_{\mu}^- \otimes u_\mu^-,0)]$
is the unique $B \times B$-invariant point of
$Z_\mu$. In particular, $Z_\mu$ contains a 
unique closed $G \times G$-orbit.
\end{Lemma}
\begin{proof}
By Borel's fixed point theorem the set of
$B \times B$-invariant points in $Z_\mu$
is nonempty. So consider a $B \times B$-invariant
point $x=[(f,a)]$, $f \in {\rm End}({\rm H}(\mu)), 
a \in k$, of $Z_\mu$. Then $f$ is a 
$B \times B$-semiinvariant element of 
${\rm End}({\rm H}(\mu)) = {\rm H}(\mu) 
\otimes   {\rm H}(\mu)^*$. In particular,
$f$ is a multiple of $v_{\mu}^- 
\otimes u_\mu^-$. As 
$v_{\mu}^- \otimes u_\mu^-$ is a weight 
vector of nonzero weight this leaves us
with two cases; either $x = [(v_{\mu}^- 
\otimes u_\mu^-,0)]$ or $x = [(0,1)]$.
Now the following Lemma \ref{closure} 
excludes the latter case and ends the 
proof of the first assertion. The second
assertion now follows as every closed 
orbit contains a $B \times B$-invariant point.
\end{proof}

\begin{Lemma}
\label{closure}
Let $M$ denote  a finite dimensional $G$-module. 
The closure in $\mathbb{P}({\rm End}(M) \oplus k)$ 
of the $G \times G$-orbit through the point 
$[(I_M,1)]$ does 
not contain the point $[(0, 1)]$.
\end{Lemma}
\begin{proof}
Let $x$ denote the point $[(I_M,1)]$. Then 
by Lem.6.1.4 in \cite{BriKum} 
$$\overline{(G \times G) \cdot x} = 
(G \times G) \overline{(T \times T) \cdot
x}.$$
In particular, it is enough to prove that 
$[(0, 1)]$ is not contained in 
$\overline{(T \times T) \cdot x} = 
\overline{T \cdot x}$.
Choose a basis $m_i$, $i \in I$, of 
$M$ consisting of $T$-eigenvectors
and let $m_i^*$ denote the dual basis.
Then, $m_i \otimes m_j^*$, $i,j 
\in I$, is a basis of ${\rm End}(M)$.
Now adjoin  $z = (0,1) \in {\rm End}(M) \oplus 
k$ to obtain a basis for  ${\rm End}(M) \oplus 
k$ and let  $z^*, m_i^* \otimes m_j$, $i,j 
\in I$, denote the dual basis. The top wedge product 
of $M$ is a 1-dimensional $G$-representation 
and hence it must be the trivial representation. 
In particular, the homogeneous polynomial function 
$$ (\prod_{i \in I} m_i^* \otimes m_i) - (z^*)^{ {\rm dim}_k(N)}$$
vanishes on $\overline{T \cdot x}$. Hence,
$[(0,1)]$ is not contained in $\overline{T \cdot x}$.
\end{proof}

\subsection{Morphisms into $Z_\mu$}

In order to state the next result we define 
$p : {\rm End}({\rm H}(\mu))^* \oplus k
\rightarrow {\rm End}({\rm H}(\mu))^*$ 
to denote the projection onto the
first summand. Recall that by Frobenius
reciprocity there exists a unique (up
to scalar) $G$-equivariant morphism  
${\rm H}(\mu)^*  \rightarrow {\rm H}(-w_0 \mu)$.

\begin{Lemma}
\label{eta}
Let $\eta : G/B \times G/B \rightarrow Z_\mu \subseteq
\mathbb P_\mu$ denote a $G \times G$-equivariant morphism.
Then 
$$\eta^*(\mathcal O_\mu(1)) \simeq 
{\mathcal L}_{G/B}(-w_0 \mu) \boxtimes 
{\mathcal L}_{G/B}(\mu)$$
and the associated map of global sections ${\rm R}_\mu$ 
fits into a commutative diagram
$$\xymatrix{ {\rm End}({\rm H}(\mu))^* \oplus k \ar[d]_p \ar[r]^{{\rm R}_\mu}& {\rm H}(-w_0 \mu)
\otimes  {\rm H}(\mu) \\
{\rm End}({\rm H}(\mu))^* \ar[r]^{\simeq} & {\rm H}(\mu)^* \otimes {\rm H}(\mu)
\ar[u] 
},$$
where the right vertical arrow is defined from a 
(unique up to nonzero scalars) nonzero 
map ${\rm H}(\mu)^*  \rightarrow {\rm H}(-w_0 \mu)$
of $G$-modules. 
\end{Lemma}
\begin{proof}
By Lemma \ref{fixed} the morphism $\eta$ must 
be given by
$$ \eta : G/B \times G/B \rightarrow 
Z_{\mu} \subseteq {\mathbb P}({\rm End}({\rm H}(\mu))
\oplus k),$$
$$ (gB, g'B) \mapsto [(gv_{\mu}^- \otimes g'u_\mu^-,0)]. $$
Hence we may, via the Segre embedding, factor $\eta$  
through the map 
$$ G/B \times G/B \rightarrow 
{\mathbb P}({\rm H}(\mu))
\times {\mathbb P}({\rm H}(\mu)^*),$$
$$ (gB, g'B) \mapsto ([(g \cdot v_\mu^-],[g' \cdot u_\mu^-]).$$
As a consequence, we find that 
$$ \eta^*({\mathcal O}_\mu(1)) \simeq 
{\mathcal L}_{G/B}(-w_0 \mu) \boxtimes 
{\mathcal L}_{G/B}(\mu).$$
Moreover the associated map on global sections 
$$ {\rm R}_\mu :{\rm End}({\rm H}(\mu))^* \oplus k = {\mathcal O}_\mu(1)({\mathbb P}_\mu) 
\rightarrow  {\rm H}(-w_0 \mu) \otimes  {\rm H}(\mu)$$
is nonzero and factorizes through the projection map $p$.
\end{proof}

\begin{Lemma}
\label{R}
The element $R_\mu (u_\mu^+ \otimes v_\mu^-,0)$
is a nonzero multiple of $v_{-w_0 \mu}^+ 
\otimes v_\mu^-$. Moreover, the element $R_\mu ({\rm Tr}_\mu ,0)$ 
is a nonzero generator of the unique diagonal $G$-invariant
line in ${\rm H}(-w_0 \mu) \otimes  {\rm H}(\mu)$.
\end{Lemma}
\begin{proof}
The first statement follows from Lemma \ref{eta}.
It follows by Frobenius reciprocity that  ${\rm H}(-w_0 \mu) 
\otimes  {\rm H}(\mu)$ contains a unique $\Delta G$-invariant
line. Hence, it suffices to show that $R_\mu ({\rm Tr}_\mu ,0)$ 
is nonzero. Consider the decomposition 
$${\rm Tr}_\mu = \sum_{\theta \in X^*(T) \times X^*(T)} 
({\rm Tr}_\mu)_\theta$$
of $Tr_\mu$ into $T \times T$-semiinvariant elements.
Then $({\rm Tr}_{\mu})_{(-w_0\mu, w_0 \mu)}$ is a nonzero 
multiple of $u_\mu^+ \otimes v_\mu^-$, and the
second statement now follows from the first 
statement.
\end{proof}

\section{Line bundles on equivariant embeddings}
  
Throughout this section we assume that $G$ is simply
connected. We use the notation introduced in Section
\ref{psi}. Let $X$ denote an equivariant embedding
of $G$.

\begin{Lemma}
\label{noboundary}
Assume that $\psi_\mu$ extends to a map 
$$ \psi_\mu : X \rightarrow \mathbb P_\mu,$$
and let $\tau_\mu$ denote the pull back 
of the global section $(u_\mu^- \otimes 
v_\mu^-,0)$ of ${\mathcal O}_\mu(1)$. 
Then the support of the divisor of zeroes 
of $\tau_\mu$ does not contain any of the 
irreducible components of $X \setminus G$.  
\end{Lemma}
\begin{proof}
Notice first that the image $\psi_\mu(X)$ 
is contained in $Z_\mu$. If the support of 
the divisor of zeroes of $\tau_\mu$ contains 
an irreducible component of  $X \setminus G$ 
then this support also contains a 
$G \times G$-orbit. Hence, also the support 
of the zero divisor of  $(u_\mu^- 
\otimes v_\mu^-,0)$ would contain a 
$G \times G$-orbit within $Z_\mu$. 
But the latter support is closed and
hence it will contain the unique closed
$G \times G$-orbit of $Z_\mu$ (see
Lemma \ref{fixed}). In particular, 
$R_\mu (u_\mu^- \otimes v_\mu^-,0)$,
and hence also $R_\mu (u_\mu^+ \otimes v_\mu^-,0)$, 
is zero.
By Lemma \ref{R} this is a contradiction.
\end{proof}

Let $D_i$ denote the closure of the Bruhat cell 
$B s_i w_0 B$ within $X$. Then 

\begin{Proposition}
\label{extension}
Assume that $\psi_{\omega_i}$ extends to a morphism 
$$ \psi_{\omega_i} : X \rightarrow \mathbb 
P_{\omega_i}.$$
Then $D_{i}$ is a locally principal Weil divisor and 
its associated line bundle ${\mathcal O}_X(D_{i})$
(resp. canonical section) is isomorphic to 
$\psi_{\omega_i}^*({\mathcal O}_{\omega_i}(1) )$
(resp. $\psi_{\omega_i}^*(u_{\omega_i}^- \otimes v_{\omega_i}^-,0)$).
\end{Proposition}
\begin{proof}
By Lemma \ref{OnG} we know that 
$\psi_{{\omega_i},G}^*(u_{\omega_i}^- \otimes v_{\omega_i}^-,0)$
is a generator of the lowest weight line in 
${\rm H}(-w_0 \omega_{i})$. Let $\hat i$ denote 
the integer satisfying $\omega_{\hat i} = - w_0
\omega_i$. Then $w_0 s_{\hat i}   =  s_i
w_0 $ and the statement is hence 
an immediate consequence of Lemma \ref{lowestfunction}
and Lemma \ref{noboundary}.
\end{proof}

As $G$ is simply connected and as $X$ is normal any 
line bundle on $X$ will have a unique $G \times G$-linearization.
This explains what is meant by an invariant global 
section in the following statement.

\begin{Corollary}
\label{divisor}
Assume that $D_{i}$ is a locally principal Weil 
divisor on $X$. Then there exists a $G \times G$-equivariant
morphism :
$$ \psi_ i : {\rm End}({\rm H}(\omega_i))^* \oplus k 
\rightarrow {\rm H}^0(X, {\mathcal O}_X(D_{i})),
$$
making the following diagram commutative
$$\xymatrix{
{\rm End}({\rm H}(\omega_i))^* \oplus k \ar[d]_{\psi_{\omega_i,G}^*}
\ar[r]^{\psi_i} & {\rm H}^0(X, {\mathcal O}_X(D_{i})) 
\ar[d]^{{\rm res}^X_G(D_i)} \\
k[G] \ar[r]^(.4){\simeq} & {\rm H}^0(G, {\mathcal O}_X(D_{i})) \\
}
$$
where the lower map is some (unique up to nonzero scalar) 
$G \times G$-equivariant identification. In particular,
there exists a nonzero $G \times G$-invariant global
section of ${\mathcal O}_X(D_{i})$.
\end{Corollary}
\begin{proof}
Consider the closure $\Gamma$ of the graph of 
$\psi_{\omega_i}$ within $X \times 
{\mathbb P}_{\omega_i}$. Projection on the 
first coordinate defines a birational projective 
morphism $\phi : \Gamma \rightarrow X$,
which is an isomorphism outside a closed 
subset of $X$ of codimension $\geq 2$
(see e.g. Prop.III.9.1,\cite{Mumford}). 
In particular, there exists a $G \times G$-stable open 
subset $X'$ of $X$ with $X \setminus X'$ 
of codimension $\geq 2$ and an extension 
of $\psi_{\omega_i}$ to $X'$. 

As $X$ is normal  
$ {\rm H}^0(X, {\mathcal O}_X(D_{i}) )
= {\rm H}^0(X', {\mathcal O}_X(D_{i})),$
which means that it suffices to prove the 
statement for $X'$. As  $\psi_{\omega_i}$ 
extends to $X'$ the first assertion now follows 
from Proposition \ref{extension}.  
To prove the second assertion it suffices to 
prove that $\psi_i(0,1)$ is nonzero. But this 
follows as  $\psi_{\omega_i,G}^*(0,1)$ is a 
nonzero (constant) function on $G$. 
\end{proof}

\section{Frobenius splittings of simply 
connected groups}

In this section $G$ denotes a simply connected 
group. The aim of this section is to construct 
a class of Frobenius splittings of $G$. These 
will be constructed by restricting Frobenius 
splittings of a certain equivariant embeddings 
$X$ of $G$. We begin by fixing the required 
properties of $X$.

\begin{Lemma}
\label{special}
There exists a smooth complete toroidal embedding
$X$ of $G$ such that the morphisms 
$\psi_{\omega_i}$, $i=1,\dots,l$,  all extend to $X$. 
\end{Lemma}
\begin{proof}
Start by choosing a complete toroidal embedding
$X'$ of $G$. Consider the product map
$\psi := \prod_{i=1}^l \psi_{\omega_i} : 
G \rightarrow \prod_{i=1}^l {\mathbb P}_{\omega_i},$
and let $\Gamma$ denote the normalization of the 
closure of the graph of $\psi$  within
$X' \times \prod_{i=1}^l {\mathbb P}_{\omega_i}$.
Then any projective resolution (see
Theorem \ref{resolution}) of $\Gamma$ has the desired 
properties. 
\end{proof} 

\begin{Remark}
A closer study of toroidal embeddings reveals 
that $\psi_{\omega_i}$ extends to arbitrary 
toroidal embedding. In particular, 
on a toroidal embedding of $G$ 
the Weil divisors $D_i$ are all locally 
principal.
\end{Remark}
 
For the rest of this section we will fix a 
toroidal embedding $X$ of $G$ satisfying
the requirements in Lemma \ref{special}
Fix a closed $G \times G$-orbit $Y$
within $X$. As noted in Section \ref{toi-emb},
we may  $G \times G$-equivariantly identify 
$Y$ with the variety $G/B \times G/B$.

\subsection{The canonical bundle}

Let $X_j, j=1, \dots, n$ denote the boundary
components of $X$, i.e. the irreducible 
components of $X \setminus G$, all of 
codimension 1 as $G$ is affine (Chap.III,\cite{Hartshorne}).  
Let furthermore $D_i$ denote the closure
of the Bruhat cell $B s_i w_0  B$ within 
$X$. By Prop.6.2.6 of \cite{BriKum} it follows that
the dual canonical bundle of $X$ is
$$ \omega_X^{-1} \simeq \mathcal{O}_X(
\sum_{j=1}^n X_j + 2\sum_{i=1}^l D_i).$$
As we will see below, the restriction 
of the line bundle  $\mathcal{O}_X(
2\sum_{i=1}^l D_i)$ to $Y$ is isomorphic
to 
${\mathcal L}_{G/B}(2 \rho) 
\boxtimes  {\mathcal L}_{G/B}(2 \rho),$
which is  the dual canonical bundle of $Y$.
Let now $\sigma_j$ denote a canonical section 
of the line bundle ${\mathcal O}_X(X_j)$ and 
let $D$ denote the divisor  $\sum_{i=1}^l 
D_i$. Then
the following statement describes the usual
way of constructing Frobenius splittings 
of $X$.

\begin{Theorem}
\label{classical}
Let $\tau$ denote a global section  
of $\mathcal{O}_X(2(p-1)D)$
such that its restriction $\tau_{|Y}$
to $Y$ corresponds to a Frobenius splitting of
$Y$. Then the global section 
$\tau \prod_{j=1}^n \sigma_j^{p-1}$
of $\omega_X^{1-p}$ corresponds (up
to a nonzero scalar) to a Frobenius
splitting of $X$.
\end{Theorem}
\begin{proof}
See e.g. proof of Thm.6.2.7 in \cite{BriKum}.
\end{proof}

\subsection{Frobenius splittings of $X$}
By the assumptions on $X$ there exists a 
commutative diagram of $G \times G$-equivariant
morphisms 
$$\xymatrix{ Y \ar[r]^(.3){\simeq} \ar@{^{(}->}[d] &  
G/B \times G/B \ar[d]^\eta \\
X \ar[r]^{\psi_{\omega_i} } & Z_{\omega_i} \\
}$$
Hence, combining Proposition \ref{extension} and 
Lemma \ref{eta} we obtain

\begin{Proposition}
\label{omega}
There exists a commutative diagram 
$$\xymatrix{ 
{\rm End}({\rm H}(\omega_i))^* \oplus 
k \ar[d]_{R_{\omega_i}} \ar[r]^{\psi_i} & {\rm H}^0(X, {\mathcal O}_X(D_{i})) 
\ar[d]^{{\rm res}_Y^X(D_{i} )} \\
{\rm H}(-w_0 \omega_i)^* \otimes {\rm H}(\omega_i) 
\ar[r]^{\simeq}  & {\rm H}^0(Y, {\mathcal O}_X(D_{i})) \\
}
$$
In particular, the image of the restriction map $res_Y^X(D_{i})$ 
contains the unique $B \times B$-stable line of 
$ {\rm H}^0(Y, {\mathcal O}_X(D_{i}))$.
\end{Proposition}

It follows from Lemma \ref{eta} that  the restriction 
of ${\mathcal O}_X((p-1)D)$ to $Y$ is isomorphic to 
the line bundle 
$ {\mathcal L}_{G/B}((p-1) \rho) \boxtimes 
{\mathcal L}_{G/B}((p-1)\rho).$  This leads 
to the following statement.

\begin{Corollary}
There exists a commutative diagram 
$$\xymatrix{ 
\bigotimes_{i=1}^l ({\rm End}({\rm H}(\omega_i))^* \oplus 
k)^{\otimes (p-1)} \ar[d]^{\otimes_i R_{\omega_i}^{p-1}} 
\ar[r]^(.6){\prod \psi_i^{p-1}} \ar[dr]^{\rm res}& 
{\rm H}^0(X, {\mathcal O}_X((p-1)D)) 
\ar[d]^{{\rm res}_Y^X((p-1)D)} \\
\bigotimes_{i=1}^l ({\rm H}(-w_0 \omega_i)^* \otimes 
{\rm H}(\omega_i))^{\otimes (p-1)} \ar[r]  &
{\rm H}^0(Y, {\mathcal O}_X((p-1) D))\\
}
$$
where the map $\rm res$ is surjective.
\end{Corollary}
\begin{proof}
By Proposition \ref{omega} the image of ${\rm res}$ is
nonzero. Hence the statement follows as 
${\rm H}^0(Y, {\mathcal O}_X((p-1) D)) \simeq {\rm St} 
\otimes {\rm St}$ is a simple $G \times G$-module 
\end{proof}

Let $v_-$ (resp. $v_+$) denote a generator of the 
lowest (resp. highest) weight space in ${\rm St}$
and let $v_\Delta$ denote a generator of the 
unique diagonal $G$-invariant line in ${\rm St}
\otimes {\rm St}$. Applying Lemma \ref{R} we find 

\begin{Lemma}
With $a_1, \dots, a_l \in k$ we have
\begin{enumerate}
\item ${\rm res}(\bigotimes_{i=1}^l (u_{\omega_i}^+ \otimes v_{\omega_i}^-,0)^
{\otimes (p-1)})$ is nonzero multiple of $v_+ \otimes v_-$.
\item ${\rm res}(\bigotimes_{i=1}^l ({\rm Tr}_{\omega i},a_i)^{\otimes (p-1)})$ is 
a nonzero  multiple of $v_{\Delta}$.
\end{enumerate} 
\end{Lemma}

By Theorem \ref{classical} and Corollary \ref{splitdiag}
we can now state and prove.

\begin{Theorem}
\label{thm-toroidal}
Let $a_1, \dots, a_l \in k$. Then the global section
$$ \prod_{i=1}^l \psi_i({\rm Tr}_{\omega i},a_i)^{p-1}
 \prod_{i=1}^l \psi_i(u_{\omega_i}^+ \otimes v_{\omega_i}^-,0) ^{p-1}
\prod_{j=1}^n \sigma_j^{p-1},$$
of $\omega_X^{1-p}$ defines a Frobenius splitting of $X$ 
(up to a nonzero constant).
\end{Theorem}

\subsection{Frobenius splittings of $G$}

The canonical bundle $\omega_G$ of $G$ is trivial and 
we may therefore choose a {\it volume form} $dG$
freely generating the global sections of $\omega_G$ as a $k[G]$-module.
Restricting the statement in Theorem \ref{thm-toroidal}
to the open subset $G$ then implies

\begin{Corollary}
\label{FrobG}
Let $a_1, \dots, a_l \in k$. Then the global section
$$ \prod_{i=1}^l \psi_{\omega_i,G}^*({\rm Tr}_{\omega i},a_i)^{p-1}
 \prod_{i=1}^l \psi_{\omega_i,G}^*(u_{\omega_i}^+ \otimes 
v_{\omega_i}^-,0)^{p-1} dG^{1-p},$$
of $\omega_G^{1-p}$ defines a Frobenius splitting of $X$ 
(up to a nonzero constant).

\end{Corollary}

\section{Frobenius splitting in the general case}

In this section $G$ denotes a simply connected linear
algebraic group. Let $X$ denote an equivariant embedding
of $G$ and let $D_i$ (resp. $\tilde{D}_i$) denote the 
closure of the Bruhat cell $B s_i w_0 B$ (resp.
$B^+ s_i B$) within X. As in the previous section
we let $X_j, j=1, \dots, n$, denote the irreducible 
components of $X \setminus G$. Moreover when 
${\bf a} \in \mathbb A^l$, we let
$\overline{\chi}^{-1}({\bf a})$ denote the closure 
of the associated Steinberg fibre within $X$.

\subsection{The smooth case} 
Assume first that $X$ is smooth. In this case 
the divisors $X_j$ are all locally principal 
and we choose canonical sections $\sigma_j$
of the associated line bundles ${\mathcal O}_X(X_i)$.
By Prop.6.2.6. in \cite{BriKum} the dual canonical bundle of
$X$ is isomorphic to 
$$ \omega_X^{-1} \simeq \mathcal{O}_X(
\sum_{j=1}^n X_j + 2\sum_{i=1}^l D_i).$$
In order to construct Frobenius splittings of
$X$ we have to construct global sections of 
$\omega_X^{1-p±}$. First we apply Corollary 
\ref{divisor} and obtain for each $i=1, \dots,
l$, a commutative diagram of 
$G \times G$-equivariant maps : 

$$\xymatrix{ 
{\rm End}({\rm H}(\omega_i))^* \oplus 
k  \ar[d]_{\psi_{\omega_i,G}^*}
\ar[r]^(.5){\psi_i} & 
{\rm H}^0(X, {\mathcal O}_X(D_i)) 
\ar[d]^{{\rm res}_G^X(D_i)} \\
k[G] \ar[r]^(.35){\simeq}  &
{\rm H}^0(G, {\mathcal O}_X(D_i))\\
}
$$
Now we have the following generalization of 
Theorem \ref{thm-toroidal} 

\begin{Corollary}
\label{cor-smooth}
Let ${\bf a} =(a_1, \dots, a_l) \in \mathbb A^l$. Then 
the global section
$$ \tau_{\bf a} = \prod_{i=1}^l \psi_i({\rm Tr}_{\omega i},-a_i)^{p-1}
 \prod_{i=1}^l \psi_i(u_{\omega_i}^+ \otimes v_{\omega_i}^-,0) ^{p-1}
\prod_{j=1}^n \sigma_j^{p-1}$$
of $\omega_X^{1-p}$ defines a Frobenius splitting of $X$ 
(up to a nonzero constant) which compatibly splits 
the closed subvarieties $\overline{\chi}^{-1}({\bf a})$,
$\tilde D_i$, $i=1,\dots, l,$ and 
$X_j$, $j=1,\dots,n$.  
\end{Corollary}
\begin{proof}
By Corollary \ref{FrobG} the restriction of $\tau_{\bf a}$
to $G$ defines a Frobenius splitting of $G$ and hence, by 
Remark \ref{opensubset}, $\tau_{\bf a}$ defines a
Frobenius splitting of $X$.  By Lemma \ref{Frob9} and 
Lemma \ref{Frob4} each component of the zero divisor
of $\tau_{\bf a}$ is compatibly Frobenius split. In
particular, the boundary divisors $X_j$ and (by 
Proposition \ref{extension}) the subvarieties 
$\tilde D_i$ are all compatibly Frobenius
split. Finally, by Lemma \ref{Frob4}, every component 
of the (scheme theoretic) intersection of the zero 
divisors of $\psi_i({\rm Tr}_{\omega i},-a_i)$
will be Frobenius split. But, by Lemma \ref{OnG}, 
the component of this latter intersection which 
intersects $G$ nontrivially is exactly the subvariety 
$\overline{\chi}^{-1}({\bf a})$.
\end{proof}

The following connected result will also be useful.

\begin{Proposition}
\label{stable}
Let ${\bf a} =(a_1, \dots, a_l) \in \mathbb A^l$. Then 
there exists a stable Frobenius splitting of $X$ along
the divisor 
$$(p-1)(\sum_{j=1}^n X_j + \sum_{i=1}^l \tilde D_i),$$
of degree 1 which compatibly Frobenius splits the 
subvariety $\overline{\chi}^{-1}({\bf a})$.
\end{Proposition}
\begin{proof}
By Lemma \ref{Frob9} and Lemma \ref{Frob5} the
Frobenius splitting $\tau_{\bf a}$ in Corollary
\ref{cor-smooth} defines a degree 1 stable Frobenius 
splitting $s$ of $X$ along the divisor 
$$(p-1)(\sum_{j=1}^n X_j + \sum_{i=1}^l \tilde 
D_i),$$ 
which compatibly Frobenius splits 
the zero divisor of the global section 
$$\prod_{i=1}^l \psi_i({\rm Tr}_{\omega i},-a_i),$$
of the line bundle $\mathcal{O}_X(\sum_{i=1}^l D_i).$
Denote by $s'$ the (by $s$) induced degree 1 stable 
Frobenius splitting of $X$ along the zero divisor $0$.
Combining  Lemma \ref{Frob2} and Lemma \ref{Frob4} 
and arguing as in the last part of the proof of
Corollary \ref{cor-smooth}, we find that the closed
subvariety $\overline{\chi}^{-1}({\bf a})$ is compatibly
Frobenius split by $s'$. By Lemma \ref{Frob1.5} it
remains to prove that $\overline{\chi}^{-1}({\bf a})$
is not contained in the support of the divisor 
$\sum_{j=1}^n X_j + \sum_{i=1}^l \tilde D_i$.
This follows from Lemma \ref{B-orbits} below.
\end{proof}

\begin{Lemma}
\label{B-orbits}
The closure of the subset $B^+ s_i  B$ within $G$
does not contain any $B$-conjugacy classes. 
\end{Lemma}
\begin{proof}
Let $f$ denote a generator of the highest weight 
line in ${\rm H}(\omega_i)$. By Lemma 
\ref{lowestfunction} the zero set of $f$ coincides
with the closure of $B^+ s_i  B$. So assume
that there exists an element $g \in G$ such that
$ f(b g b^{-1}) = 0$ for all $b \in B$. By 
definition of ${\rm H}(\omega_i)$ this 
implies that $f(b g)= 0$ for all 
$b \in B$. Moreover, as $f$ is $B^+$-semiinvariant 
we find that $f(h g) = 0$ for 
every element $h$ in the open dense subset 
$B^+ B$ of $G$. In particular, $f$ vanishes
on an open dense subset of $G$ which is 
a contradiction.
\end{proof}

\begin{Corollary}
\label{stable-ample}
Assume that $X$ is projective and let  
${\bf a} =(a_1, \dots, a_l) \in \mathbb A^l$.  
Then there exists a stable Frobenius splitting of $X$ 
along an ample divisor with support $X \setminus G$ 
which compatibly Frobenius 
splits the subvariety $\overline{\chi}^{-1}({\bf a})$.
\end{Corollary}
\begin{proof}
By Proposition \ref{stable} and Lemma \ref{Frob2} 
there exists a stable Frobenius splitting of $X$ 
along the divisor $\sum_{j=1}^n X_j$ which 
compatibly splits $\overline{\chi}^{-1}({\bf a})$.
Applying  Lemma \ref{Frob2} and Lemma \ref{Frob3} 
it suffices to show that there exists positive 
integers $c_j>0$ such that $\sum_{j=1}^n c_j X_j$ 
is ample. This follows from Prop.4.1(2) in 
\cite{BriTho}.
\end{proof}

\subsection{The general case}

For a general equivariant $G$-embedding we can
now prove.

\begin{Theorem}
\label{general-split}
Let $X$ denote an arbitrary equivariant 
$G$-embedding and let ${\bf a} =(a_1, \dots, a_l) \in \mathbb A^l$. 
Then $X$ is Frobenius splits compatibly with the 
closed subvarieties $\overline{\chi}^{-1}({\bf a})$,
$\tilde D_i$, $i=1,\dots, l$ and $X_j$, $j=1,\dots,n$.  
\end{Theorem}
\begin{proof}
By Theorem \ref{resolution} there exists a 
resolution $f : X' \rightarrow X$ of $X$ 
by a smooth $G$-embedding $X'$. By Zariski's
main theorem we know
$f_* {\mathcal O}_{X'} \simeq {\mathcal O}_{X}$. 
Hence, by the discussion in Section 
\ref{push-forward} it suffices to prove 
the above statement for $X'$. Now apply 
Corollary \ref{cor-smooth}.
\end{proof}

\subsection{Canonical Frobenius splittings}

Let $Z$ denote an arbitrary $G$-variety and let 
$s : F_* \mathcal O_Z \rightarrow \mathcal O_Z$
denote a Frobenius splitting of $Z$. For any 
root $\alpha$ we let $x_\alpha : k \rightarrow
G$ denote the associated root homomorphism
satisfying $t \in T, c \in k : t x_\alpha(c)
t^{-1} = x_\alpha(\alpha(t) c)$. Recall
(Defn.4.1.1,\cite{BriKum}) that $s$ is 
said to be a {\it canonical Frobenius splitting} 
if $s$ is $T$-invariant and satisfies
$$ x_{\alpha_i}(c) s = \sum_{j=1}^{p-1} 
c^j s_j$$
for every simple root $\alpha_i$, and certain
elements $s_j \in {\rm Hom}_{\mathcal O_X} (F_* \mathcal O_Z,
\mathcal O_Z)$. The primary reason for the interest in 
canonical Frobenius splittings comes from 
the following consequence (see e.g. Thm.4.2.13
in \cite{BriKum}) 

\begin{Theorem}
Assume that $Z$ has a canonical Frobenius splitting.
Let $\mathcal L$ denote a  $G$-linearized line 
bundle on $Z$. Then the $G$-module ${\rm H}^0(Z, \mathcal L)$ 
admits a good filtration, i.e. there exists 
a sequence of $G$-modules
$$0 =M^0 \subseteq M^1 \subseteq M^2 \subseteq \cdots $$
such that ${\rm H}^0(Z, \mathcal L) = \cup_i M^i$
and satisfying that the successive quotients $M^{j+1}/M^{j}$ are 
isomorphic to modules of the form ${\rm H}(\lambda_j)$
for certain dominant weights $\lambda_j$.
\end{Theorem}

Let ${\bf a} =(a_1, \dots, a_l) \in \mathbb A^l$ and
consider an arbitrary equivariant embedding $X$ of 
$G$. Then there exists a diagonal $G$-action on the 
closed subvariety $\overline{\chi}^{-1}({\bf a})$.
We claim

\begin{Corollary}
The variety $\overline{\chi}^{-1}({\bf a})$ admits a
canonical Frobenius splitting. In particular,
when $\mathcal L$ is a $G$-linearized line 
bundle on  $\overline{\chi}^{-1}({\bf a})$ then 
the $G$-module  ${\rm H}^0(\overline{\chi}^{-1}({\bf a}) 
, \mathcal L)$ has a good filtration.
\end{Corollary}
\begin{proof}
By the results in the previous section 
$\overline{\chi}^{-1}({\bf a})$ is compatibly
Frobenius split by a Frobenius splitting 
$s$ of $X$. Hence, it suffices to prove 
that $s$ is a canonical Frobenius splitting
of $X$. By the proof of Theorem 
\ref{general-split} it moreover suffices to 
consider the case when $X$ is smooth and 
$s$ is defined from $\tau_{\bf a}$ 
(with notation as in Corollary \ref{cor-smooth}).
Now every factor of $\tau_{\bf a}$ except 
$$\prod_{i=1}^l \psi_i(u_{\omega_i}^+ 
\otimes v_{\omega_i}^-,0) ^{p-1}$$
is invariant under the diagonal $G$-action.
Hence we may concentrate on the (diagonal) 
$T$-invariant factor $\psi_i(u_{\omega_i}^+ 
\otimes v_{\omega_i}^-,0)$. The statement
now follows as, for all $j$, 
$$ x_{\alpha_j}(c) v_{\omega_i}^- = v_{\omega_i}^- + 
c v_{i,j},$$   
$$ x_{\alpha_j}(c) u_{\omega_i}^+ = u_{\omega_i}^+,$$ 
for certain elements $v_{i,j} \in {\rm H}(\omega_i)$
(recall that the $(w_0 \omega_i + q \alpha_j)$-weight
space of ${\rm H}(\omega_i)$ is zero for $q >1$). 
\end{proof}

\section{Cohomology vanishing}

In this section we will discuss various results 
which will enable us to prove that closures of 
Steinberg fibers in arbitrary equivariant 
embeddings have nice geometric properties

\subsection{Resolutions and direct images}

Let $X'$  denote a projective equivariant embedding 
of $G$ and let $f : X \rightarrow X'$  be a resolution 
of $X'$ by a smooth projective equivariant embedding.
When ${\bf a} = (a_1, \dots, a_l) \in \mathbb A^l$ we denote by 
$\overline{\chi}^{-1}({\bf a})$ (resp. 
$\overline{\chi}^{-1}({\bf a})'$) the closure 
of the Steinberg fiber at ${\bf a}$ within
$X$ (resp. $X'$).

\begin{Corollary}
\label{trivial}
With the notation described above we have 
\begin{enumerate}
\item $f_* {\mathcal O}_X = {\mathcal O}_{X'}$ and 
$R^i f_* {\mathcal O}_X = 0$ when $i>0$.
\item $f_* {\mathcal O}_{ \overline{\chi}^{-1}({\bf a})} 
= {\mathcal O}_{\overline{\chi}^{-1}({\bf a})'}$ 
and $R^i f_* {\mathcal O}_ {\overline{\chi}^{-1}({\bf a})}= 0$ when $i>0$.
\end{enumerate}
\end{Corollary}
\begin{proof}
As $X'$ is normal and $f$ is birational it 
follows from Zariski's main theorem that 
$f_* \mathcal O_X = \mathcal O_{X'}$. 
Hence, by Lemma \ref{Frob8} it suffices
to prove that there exists a stable 
Frobenius splitting of $X$ along an 
ample divisor which compatibly Frobenius 
splits $\overline{\chi}^{-1}({\bf a})$.
Now apply  Corollary \ref{stable-ample}.
\end{proof}

\subsection{Cohomology}

We are now ready to prove the following statement 
about cohomology of line bundles on closures of 
Steinberg fibers. 

\begin{Proposition}
\label{vanishing}
Let $X$ denote a projective equivariant embedding
of $G$ and let ${\bf a} = (a_1,\dots, a_l) \in 
\mathbb A^l$. Let 
$\mathcal M$ (resp. $\mathcal L$) denote a globally
generated line bundle on $X$ (resp. 
$\overline{\chi}^{-1}({\bf a})$). Then
$${\rm H}^i(X, {\mathcal M})= {\rm H}^i(\overline{\chi}^{-1}
({\bf a}), {\mathcal L}) = 0~,i>0.$$ 
Moreover, the restriction map 
$${\rm H}^0(X, {\mathcal L}) \rightarrow  
{\rm H}^0(\overline{\chi}^{-1}({\bf a}) 
, {\mathcal L}),$$ 
is surjective.
\end{Proposition}
\begin{proof}
By Corollary \ref{trivial} we may assume that 
$X$ is smooth. Now apply Corollary \ref{stable-ample}
and the ``in particular'' parts of Lemma 
\ref{Frob6} and Lemma \ref{Frob7}.
\end{proof}
 
\section{Geometry of closures of Steinberg fibers}

In this section we will study the geometry of the
closure of a Steinberg fiber within an equivariant
$G$-embedding $X$ of a simply connected group $G$. 

\subsection{The smooth case}

Assume that $X$ is a smooth $G$-embedding and 
let ${\bf a}= (a_1,\dots, a_l) \in \mathbb A^l$. 
Let $\mathcal{O}_{\bf a}$ denote the restriction 
of the line bundle 
$$\mathcal{O}_X(\sum_{j=1}^n X_j + \sum_{i=1}^l \tilde D_i)$$
to $\overline{\chi}^{-1}({\bf a})$ and and let 
$s_{\bf a}$ denote the associated restricted canonical 
section. Notice that by Proposition \ref{stable}
the section $s_{\bf a}$ is nonzero.

\begin{Theorem}
\label{lci}
Let $X$ denote a smooth $G$-embedding and let  
${\bf a} = (a_1, \dots,a_l) \in {\mathbb A}^l$.
Then the closure $\overline{\chi}^{-1}({\bf a})$
of the Steinberg fiber $\chi^{-1}({\bf a})$ 
within $X$ is a local complete intersection. 
Moreover, the dualizing sheaf on $\overline{\chi}^{-1}
({\bf a})$ equals the line bundle 
${\mathcal O}_{\bf a}^{-1}$. 
\end{Theorem}
\begin{proof}
As $\omega_X^{-1} = \mathcal O_X(\sum_j X_j + \sum_i (D_i
+ \tilde D_i))$ it suffices to prove that $\overline{\chi}^{-1}({\bf a})$
coincides with the intersection $Z$ of
the zero schemes of the sections  
$\psi_i({\rm Tr}_{\omega i},-a_i)$,
$i=1, \dots, l$. When $X=G$ this is clearly the case. Hence, it suffices to 
prove that all components of $Z$ intersect the 
open subset $G$. So assume that $Z'$ is a component
of $Z$ which does not intersect $G$. Then $Z'$
must be contained in one of the boundary components
$X_j$. Choose an open affine subset $U$ 
such that $\emptyset \neq U \cap Z' = U \cap Z$
and such that the canonical bundle $\omega_U$
is trivial. Choose also a volume form $dU$, i.e.
a global section of $\omega_U$ generating the 
set of global sections as a $k[U]$-module. 
Then $\tau_{\bf a}$ (with notation as in Corollary 
\ref{cor-smooth}) restricts to a Frobenius splitting 
of $U$ of the form 
$$ (\tau_{\bf a})_{|U} = (\prod_{i=1}^l f_i^{p-1}
\prod_{i=1}^l g_i^{p-1} \prod_{j=1}^n h_j^{p-1}) 
\cdot (dU)^{1-p}$$
where $f_i$, $g_i$ and $h_i$ are the functions on 
$U$ defining the restrictions of the zero sets of 
$\psi_i({\rm Tr}_{\omega i},-a_i)$, $\psi_i(
u_{\omega_i}^+ \otimes v_{\omega_i}^-,0)$ and 
$\sigma_j$ respectively.
By assumption, the common zero $V(f_1,\dots,f_l)$
of $f_1, \dots, f_l$ is contained in the zero set
$V(h_j)$ of $h_j$. Hence by Hilberts Nullstellensatz, 
$h_j$ is contained in the radical of the ideal 
$(f_1,\dots,f_l)$. But $(f_1,\dots,f_l)$ is compatibly
split by Lemma \ref{Frob9} and Lemma \ref{Frob4}
and therefore it is radical by Lemma \ref{Frob4}. 
We conclude, that $h_j \in (f_1,\dots,f_l)$
and hence 
$$\prod_{i=1}^l f_i^{p-1}
\prod_{i=1}^l g_i^{p-1} \prod_{j=1}^n h_j^{p-1}
\in (f_1^p, \dots, f_l^p).$$
In particular, the Frobenius splitting $(\tau_{\bf a})_{|U}$
of $U$ will map the constant function $1$ to an
element within the ideal $(f_1,\dots, f_l)$. Therefore
$1 \in (f_1,\dots, f_l)$ and hence 
$$Z' \cap U = V(f_1,\dots, f_l) = \emptyset,$$
which is a contradiction.
\end{proof}

\begin{Corollary}
\label{smooth-geometry}
Let $X$ denote a smooth $G$-embedding and let  
${\bf a} = (a_1, \dots,a_l) \in {\mathbb A}^l$.
Then the closure $\overline{\chi}^{-1}({\bf a})$
of the Steinberg fiber $\chi^{-1} ({\bf a})$ within 
$X$ is normal, Gorenstein and Cohen-Macaulay. 
\end{Corollary}
\begin{proof}
By Theorem \ref{lci} it suffices to show that 
the closure $\overline{\chi}^{-1} ({\bf a})$
is smooth in codimension 1. So let 
$Z$ denote an irreducible component of the singular 
locus of $\overline{\chi}^{-1} ({\bf a})$. 
If $G \cap Z
\neq \emptyset$ then the codimension of $Z$ 
is $\geq 2$ as $\chi^{-1} ({\bf a})$ is normal
by Thm.6.11 in \cite{Steinberg}. So assume that 
$Z$ is contained in one of the boundary 
components $X_j$ of $X$. 

Consider the scheme theoretic intersection
$X_j \cap \overline{\chi}^{-1} ({\bf a})$
which is reduced by Corollary \ref{cor-smooth}
and Lemma \ref{Frob4}.
Hence, as $X_j$ is locally principal every 
smooth point of $X_j \cap \overline{\chi}^{-1} 
({\bf a})$ will also be a smooth point 
of $\overline{\chi}^{-1} ({\bf a})$. In 
particular, $Z$ is properly contained 
in one of the components of $X_j \cap 
\overline{\chi}^{-1} ({\bf a})$.
But the components of the latter scheme
all have codimension 1 in 
$\overline{\chi}^{-1} ({\bf a})$
and hence $Z$ must have codimension
$\geq 2$ in $\overline{\chi}^{-1} ({\bf a})$. 
This ends the proof.
\end{proof}

\subsection{General $G$-embeddings}

Now assume that $X$ is an arbitrary equivariant 
embedding
of a simply connected group $G$. The following 
result is due to G. Kempf. The version stated 
here is taken from \cite{BriPol}.

\begin{Lemma}
\label{Kempf}
Let $f : Z' \rightarrow Z$ denote a proper map 
of algebraic schemes satisfying that $f_ *
{\mathcal O}_{Z'} = {\mathcal O}_{Z}$ and 
$R^if_* {\mathcal O}_{Z'} = 0,  ~i>0$.  
If $Z'$ is Cohen-Macaulay with dualizing sheaf
$\omega_{Z'}$ and if $R^if_* {\omega}_{Z'} = 0$
for $i >0$, then $Z$ is Cohen-Macaulay with 
dualizing sheaf $f_* \omega_{Z'}$.
\end{Lemma}

We will also need the following result due
to Mehta and van der Kallen (\cite{MehKal})

\begin{Lemma}
\label{Meh-Kal}
Let $f : Z' \rightarrow Z$ denote a proper 
morphism of schemes and let $V'$ (resp. $V$) 
denote a closed subscheme of $Z'$ (resp. $Z$).
Let ${\mathcal I}_{V'}$ denote the sheaf of 
ideals of $V'$. Fix an integer $i$ and assume
\
\begin{enumerate}
\item $f^{-1} (V) \subseteq V'$.
\item $R^if_*{\mathcal I}_{V'}$ vanishes outside
$V$.
\item $V'$ is compatibly Frobenius split in $Z'$.
\end{enumerate}
Then $R^if_*{\mathcal I}_{V'}= 0$.
\end{Lemma}

We are ready to prove

\begin{Theorem}
Let $X$ denote an arbitrary equivariant embedding of 
$G$ and let ${\bf a} = (a_1, \dots, a_l) \in \mathbb
A ^l$. Then the closure $\overline{\chi}^{-1}({\bf a})$
 of the Steinberg fiber 
at ${\bf a}$ in $X$ is normal and Cohen-Macaulay.
\end{Theorem}
\begin{proof}
Any equivariant embedding has an open cover by 
open equivariant subsets of projective equivariant
embeddings. (see e.g. proof of Cor.6.2.8 in \cite{BriKum}) 
This reduces the statement to the case where
$X$ is projective. Choose a projective 
resolution $f : X' \rightarrow X$ of $X$ 
by a smooth equivariant embedding $X'$. 
Then
$f_* {\mathcal O}_{\overline{\chi}^{-1}
({\bf a})'}= {\mathcal O}_{\overline{\chi}^{-1}({\bf a})}$ 
and $R^if_* {\mathcal O}_{\overline{\chi}^{-1}
({\bf a})'} = 0,  ~i>0$, by Corollary \ref{trivial}. 
Hence by Corollary \ref{smooth-geometry} this 
implies that $\overline{\chi}^{-1}({\bf a})$ 
is normal. 

In order to show that $\overline{\chi}^{-1}({\bf a})$ 
is Cohen-Macaulay we apply the above Lemma \ref{Kempf}
and Lemma \ref{Meh-Kal}. By Theorem 
\ref{lci} it suffices to prove that $R^if_* 
{\mathcal O'}_{{\bf a}}^{ -1} = 0$, $i>0$. 
Let $Z$ denote the closed subscheme of 
$\overline{\chi}^{-1}({\bf a})'$ defined 
by $s_{\bf a}'$. Combining Proposition 
\ref{stable} and Lemma \ref{Frob5} we 
find that $\overline{\chi}^{-1}({\bf a})'$ 
is Frobenius split compatibly with the 
closed subscheme $Z$. Moreover, $f : 
\overline{\chi}^{-1}({\bf a})'
\rightarrow \overline{\chi}^{-1}({\bf a})$ is an 
isomorphism above the the open subset $\chi^{-1}({\bf a})$
and 
$$f^{-1}(\overline{\chi}^{-1}({\bf a}) \setminus 
\chi^{-1}({\bf a})) \subseteq Z$$
Hence, by Lemma \ref{Meh-Kal} we conclude
$R^if_*{\mathcal I}_{Z}= 0$ for $i>0$. But 
${\mathcal I}_{Z} \simeq {\mathcal O'}_{\bf a}^{-1}$
which ends the proof.
\end{proof}

\begin{Remark}
Let $\overline{\chi}^{-1}({\bf a})$ denote the closure 
of the Steinberg fiber at a point ${\bf a} \in \mathbb A^l$ 
within some equivariant 
embedding $X$ of $G$. Similar as to the situation for large 
Schubert varieties (see \cite{BriTho}) we expect that 
 $\overline{\chi}^{-1}({\bf a})$ is strongly $F$-regular
(see \cite{Smith}). In order to prove this
it suffices to prove that $\overline{\chi}^{-1}({\bf a})$ 
is globally $F$-regular (see \cite{Smith}) when $X$ 
is a projective smooth embedding. In this case we have 
seen that $\overline{\chi}^{-1}({\bf a})$ is compatibly 
Frobenius split by a stable Frobenius splitting of $X$ 
along an ample divisor $D$ with support $X \setminus G$ 
(Corollary \ref{stable-ample}). By Thm.3.10 in \cite{Smith} 
it therefore suffices to prove that ${\chi}^{-1}({\bf a})$ 
is strongly $F$-regular. When ${\chi}^{-1}({\bf a})$ 
coincides with the conjugacy class of a regular semisimple
element of $G$ this is clearly the case. However, for the 
most interesting fiber of  $\chi$, i.e. the unipotent
variety $\mathcal U$ of $G$, we do not know how to prove
the latter statement.
\end{Remark}

\bibliographystyle{amsplain}

\end{document}